\documentclass[12pt]{amsart}
\usepackage{pgf,tikz,pgfplots} %For TikZ
\pgfplotsset{compat=1.14}
\usepackage{mathrsfs}
\usetikzlibrary{arrows}
\usetikzlibrary{patterns}
\usepackage{pgfplots}
\usetikzlibrary{intersections, pgfplots.fillbetween}
\usepackage[colorlinks=true,urlcolor=blue, citecolor=red,linkcolor=blue,linktocpage,pdfpagelabels, bookmarksnumbered,bookmarksopen]{hyperref}
\usepackage[hyperpageref]{backref}
\usepackage{cleveref}
\usepackage{xcolor}
\usepackage{amsthm} %for citing inside theorem header
\usepackage{latexsym,amsmath,amssymb}
\usepackage{accents}
\usepackage{a4wide}
\usepackage{soul}
\usepackage{mathtools} %For \chi with a much lower index (\mychi)
\usepackage{xparse} %For a new definiton of \chi with a slightly lower index

\pgfdeclarelayer{ft}
\pgfdeclarelayer{bg}
\pgfsetlayers{bg,main,ft}

\usepackage[colorinlistoftodos,prependcaption,textsize=tiny]{todonotes}

\usepackage{tikz}

\definecolor{indigo}{rgb}{0.29, 0.0, 0.51}
\definecolor{p1}{gray}{0.4}
\definecolor{p2}{gray}{0.6}
\definecolor{p3}{gray}{0.98}
\definecolor{p4}{gray}{0.8}
\definecolor{p5}{gray}{0.9}

plus 6pt minus 12pt
\def\eps{\varepsilon}

\def\B{\mathbb{B}}

\def\S{{\mathbb S}}

\renewcommand{\div}{{\rm div}}

\newtheorem{theorem}{Theorem}
\newtheorem{lemma}[theorem]{Lemma}

\newtheorem{proposition}[theorem]{Proposition}

\newtheorem{definition}[theorem]{Definition}
\newtheorem{example}[theorem]{Example}

\newtheorem{condition}{Condition}

\crefname{condition}{condition}{conditions}
\Crefname{condition}{Condition}{Conditions}

\newcommand{\R}{\mathbb{R}}

\newcommand{\barint}{
\rule[.036in]{.12in}{.009in}\kern-.16in \displaystyle\int }

\newcommand{\barcal}{\mbox{$ \rule[.036in]{.11in}{.007in}\kern-.128in\int $}}

\def\mvint_#1{\mathchoice
          {\mathop{\vrule width 6pt height 3 pt depth -2.5pt
                  \kern -8pt \intop}\nolimits_{\kern -3pt #1}}%
          {\mathop{\vrule width 5pt height 3 pt depth -2.6pt
                  \kern -6pt \intop}\nolimits_{#1}}%
          {\mathop{\vrule width 5pt height 3 pt depth -2.6pt
                  \kern -6pt \intop}\nolimits_{#1}}%
          {\mathop{\vrule width 5pt height 3 pt depth -2.6pt
                  \kern -6pt \intop}\nolimits_{#1}}}

\numberwithin{theorem}{section} \numberwithin{equation}{section}

\newcommand{\aleq}{\precsim}
\newcommand{\ageq}{\succsim}

\let\latexchi\chi
\makeatletter
\renewcommand\chi{\@ifnextchar_\sub@chi\latexchi}
\newcommand{\sub@chi}[2]{% #1 is _, #2 is the subscript
  \@ifnextchar^{\subsup@chi{#2}}{\latexchi_{#2}}%
}
\newcommand{\subsup@chi}[3]{% #1 is the subscript, #2 is ^, #3 is the superscript
  \latexchi_{#1}^{#3}%
}
\makeatother
\makeatletter
\def\tikz@arc@opt[#1]{% over-write!
  {%
    \tikzset{every arc/.try,#1}%
    \pgfkeysgetvalue{/tikz/start angle}\tikz@s
    \pgfkeysgetvalue{/tikz/end angle}\tikz@e
    \pgfkeysgetvalue{/tikz/delta angle}\tikz@d
    \ifx\tikz@s\pgfutil@empty%
      \pgfmathsetmacro\tikz@s{\tikz@e-\tikz@d}
    \else
      \ifx\tikz@e\pgfutil@empty%
        \pgfmathsetmacro\tikz@e{\tikz@s+\tikz@d}
      \fi%
    \fi
    \tikz@arc@moveto
    \xdef\pgf@marshal{\noexpand%
    \tikz@do@arc{\tikz@s}{\tikz@e}
      {\pgfkeysvalueof{/tikz/x radius}}
      {\pgfkeysvalueof{/tikz/y radius}}}%
  }%
  \pgf@marshal%
  \tikz@arcfinal%
}
\let\tikz@arc@moveto\relax
\def\tikz@arc@movetolineto#1{%
  \def\tikz@arc@moveto{\tikz@@@parse@polar{\tikz@arc@@movetolineto#1}(\tikz@s:\pgfkeysvalueof{/tikz/x radius} and \pgfkeysvalueof{/tikz/y radius})}}
\def\tikz@arc@@movetolineto#1#2{#1{\pgfpointadd{#2}{\tikz@last@position@saved}}}
\tikzset{%
  move to start/.code=\tikz@arc@movetolineto\pgfpathmoveto,%
  line to start/.code=\tikz@arc@movetolineto\pgfpathlineto}
\makeatother
\usepackage[shortlabels]{enumitem}
\usepackage{url}
\usepackage{caption}

\author{Armin Schikorra}
\address[Armin Schikorra]{Department of Mathematics,
University of Pittsburgh,
301 Thackeray Hall,
Pittsburgh, PA 15260, USA}
\email{armin@pitt.edu}
\author[Akshara Vincent]{Akshara Vincent}
\address[Akshara Vincent]{Department of Mathematics,
 University of Pittsburgh,
 301 Thackeray Hall,
 Pittsburgh, PA 15260, USA}
 \email{akv20@pitt.edu}

\title{Calder\'on--Zygmund estimates and their failure for PDE with linear growth }
\begin{document}
\begin{abstract}
We investigate when Calder\'on--Zygmund estimates hold for PDEs of the type 
\[
{\rm div} (A(\nabla u)) = \div (F)
\]
where $A: \R^2 \to \R^2$ is a continuous map with linear growth. We show that, under natural assumptions, $L^p$-Calder\'on--Zygmund estimates hold for all $p \in (1,\infty)$ if and only if $v \mapsto \lim_{t \to \infty} A(tv)/t$ is linear and non-degenerate. 
\end{abstract}
\maketitle
\tableofcontents

\section{Introduction}
% Consider (regular) solutions $u$ to the \emph{linear} equation
% \[
% \begin{cases}
% \div(A \nabla u) = \div(F) \quad &\text{in $\B^n$},\\
% u = 0 \quad &\text{on $\partial \B^n$}.
% \end{cases}
% \]
% If $A: \B^n \to \R^{n \times n}_{sym}$ is symmetric, bounded, and uniformly elliptic and belongs to $VMO(\R^n)$ then we have Calder\'on--Zygmund estimates 
% \begin{equation}\label{eq:CZest}
% \|\nabla u\|_{L^p(\B^n)} \aleq_{A,p} \|F\|_{L^p(\B^n)} 
% \end{equation}
% hold for all $p \in (1,\infty)$, see, e.g. \cite{IS98}. If $A$ has less regularity, e.g. is only bounded measurable then Calder\'on--Zygmund estimates such as \eqref{eq:CZest} may not hold for all $p \in (1,\infty)$, but may hold only in a range $p \in (2-\eps,2+\eps)$ -- where $\eps$ depends on the ellipticity on $A$, see \cite{Meyers63,AFS08} and also \cite{Kamthorntaksina2025,MolinaBakhos2026}.

Consider (regular) solutions $u$ to the \emph{linear} equation
\[
\begin{cases}
\div(A \nabla u) = \div(F) \quad &\text{in $\B^n$},\\
u = 0 \quad &\text{on $\partial \B^n$}.
\end{cases}
\]
If $A: \B^n \to \R^{n \times n}_{sym}$ is symmetric, bounded, uniformly elliptic, and has $VMO$-regularity, then the Calder\'on--Zygmund estimate
\begin{equation}\label{eq:CZest}
\|\nabla u\|_{L^p(\B^n)} \aleq_{A,p} \|F\|_{L^p(\B^n)}
\end{equation}
holds for all $p \in (1,\infty)$; see, for example, \cite{IS98}. If $A$ has less regularity, for example is only bounded measurable, then Calder\'on--Zygmund estimates such as \eqref{eq:CZest} may not hold for all $p \in (1,\infty)$, but may hold only in a range $p \in (2-\eps,2+\eps)$, where $\eps$ depends on the ellipticity of $A$ (see \cite{AFS08, Meyers63}, and also \cite{Kamthorntaksina2025,MolinaBakhos2026}).

% It is a natural question to consider the nonlinear version, i.e. the situation when $v \mapsto A(v)$ is not linear, but still smooth and of linear growth. 
% Motivated by recent positive and negative results, \cite{BDS,colombo2022nonclassicalsolutionsplaplaceequation,Johansson24,Armin24,Vincent25}, we are not merely interested in examples for $A$ where Calder\'on--Zygmund estimate fail for some $p$ -- but rather in a precise characterization when this happens. 

It is a natural question to consider the nonlinear version, that is, the situation when $v \mapsto A(v)$ is not linear, but still smooth and of linear growth. Motivated by recent positive and negative results \cite{BDS,colombo2022nonclassicalsolutionsplaplaceequation,Johansson24,Armin24,Vincent25}, we are not merely interested in examples for $A$ where the Calder\'on--Zygmund estimate fails for some $p$, but rather in a precise characterization of when this happens.

Here is our main result.

%==============================================================
\begin{theorem}\label{th:main}
    Let $A:\R^2\to \R^2$ be a continuous map and assume there exists some $L>0$ such that
    \begin{equation}\label{eq:linearity of A}
        |A(v)|\leq L|v|\quad \forall\;v\in \R^2.
    \end{equation}
    Suppose there exists $H:\S^1\to \R^2$ such that
    \begin{equation}\label{eq:limit of A}
        \lim_{t\to \infty} \sup_{v\in \S^1}\left| \frac{A(tv)}{t}-H(v)\right|=0. 
    \end{equation}
    Then, the following are equivalent. 
    \begin{enumerate} \item \label{item 1:thm main}For every $1<p<\infty$, there exists a constant $C>0$ (depending only on $A$ and $p$) such that
    \begin{equation}\label{eq:CZ est main thm}
        \int_{\B^2}|\nabla u|^p \leq C \left(1+\int_{\B^2} |F|^p\right)
    \end{equation}
    whenever $u\in W^{1,p}_0(\B^2)$ and $F\in L^p(\B^2,\R^2)$ solve
    $$\div(A(\nabla u))=\div (F)\quad \text{in}\;\B^2.$$
    \item \label{item 2 them main} The map $H$ satisfies 
    \begin{equation}\label{eq:nondegeneracy}
        H(v)\cdot v\neq 0\quad \text{for } v\in \S^1
    \end{equation}
    and it is extendable to a linear map $\Tilde{H}:\R^2\to \R^2$ given by  
    \begin{equation}\label{eq:extnesion of H}
        \Tilde{H}(v)=\begin{cases}
            |v| H\left(\frac{v}{|v|}\right),\quad & v\neq 0\\
            0, \quad & v=0.
        \end{cases}
    \end{equation}
    \end{enumerate}
\end{theorem}

One direction of \Cref{th:main} is a consequence of the following two results, which might be interesting on their own -- the first one is very much expected.
%===============================================================
\begin{theorem}\label{thm:coercive counterexample}
    Let $A:\R^2\to \R^2$ be a continuous map satisfying \eqref{eq:linearity of A} and \eqref{eq:limit of A}. Assume that there exists a $v_0 \in \S^1$ such that $H(v_0)\cdot v_0 =0$. Then, for any $\Lambda >0$ and $p\in (1,\infty)$, there exist a Lipschitz function $u\in W^{1,\infty}_0(\B^2)$ and a map $F\in L^\infty(\B^2, \R^2)$ such that
\begin{equation}\label{eq:PDE_thmA}
    \begin{cases}
\div (A(\nabla u)) = \div F \quad& \text{in $\B^2$},\\
u=0 \quad &\text{on $\partial \B^2$.}
\end{cases}
\end{equation}
    But
    \begin{equation}\label{eq:CZ_estimate_thmA}
        \int_{\B^2}|\nabla u|^p>\Lambda\left(1+ \int_{\B^2}|F|^p\right).
    \end{equation}
\end{theorem}
%=================================================================

The second, probably more interesting, property relies on:
\begin{condition}\label{cond:A}
Let $A:\mathbb{R}^2\to\mathbb{R}^2$ be a continuous map satisfying the following: \\
There exist non-zero vectors
\[
v_1,v_2,v_3,w_1,w_2,w_3\in\R^2,
\]
and a number $\beta\in (0,1)$ such that
%%%%%%%%%%%%%%%%%%%%%%%%%%%%%%%%%%%%%%%%%%%
\begin{figure}[ht]
\centering
\begin{minipage}{0.55\textwidth}
  \begin{equation}\label{eq:unit_vector_sum_v}
v_1+v_2+v_3=0,
\end{equation}
\begin{equation}\label{eq:unit_vector_sum_w}
w_1+w_2+\beta w_3=0,
\end{equation}
\begin{equation}\label{eq:coercivity}
v_i\cdot w_i>0 \quad \text{for}\;\;i=1,2,3,
\end{equation}
and 
\begin{equation}\label{eq:three_ray_conv}
\lim_{t\to \infty}\frac{A(tv_i)}{t}= w_i
\quad \text{for}\;\;i=1,2,3.
\end{equation}
\end{minipage}
\hfill
\begin{minipage}{0.42\textwidth}
  \centering
  \begin{tikzpicture}[scale=1.5, >=stealth, thick]
    \fill (0,0) circle (0.6pt);

    \draw[black, ->] (0,0) -- (44:1)   node[above right] {$v_1$};
    \draw[black, ->] (0,0) -- (-32:1)  node[below right] {$v_2$};
    \draw[black, ->] (0,0) -- (179:1)  node[left]       {$v_3$};

    \draw[blue, ->] (0,0) -- (22:1)    node[right]       {$w_1$};
    \draw[blue, ->] (0,0) -- (-56:1)   node[below right] {$w_2$};
    \draw[blue, ->] (0,0) -- (162:1.9) node[left]  {$w_3$};
  \end{tikzpicture}
  \caption{ Vectors $v_i$ and $w_i$}
\end{minipage}

\end{figure}
%%%%%%%%%%%%%%%%%%%%%%%%%%%%%%%%%%%%%%%%%%

\end{condition}

\begin{theorem}\label{th:CZ_failure_A}
Let $A:\R^2\to\R^2$ be continuous and assume that, for some $0<L<\infty$,
\begin{equation}
 |A(\xi)|\le L |\xi|
\qquad \forall\; \xi\in\R^2.
\end{equation}
Assume that $A$ satisfies \Cref{cond:A}. Then, there exists a $\bar{q}>1$ such that the following holds: For any $\Lambda >0$ and $p\in (1,\bar{q})$, there exist a Lipschitz function $u\in W^{1,\infty}_0(\B^2)$ and a map $F\in L^\infty(\B^2, \R^2)$ satisfying \eqref{eq:PDE_thmA} and \eqref{eq:CZ_estimate_thmA}.
% \begin{equation}\label{eq:PDE_thmA}
%     \begin{cases}
% \div (A(\nabla u)) = \div F \quad& \text{in $\B^2$},\\
% u=0 \quad &\text{on $\partial \B^2$.}
% \end{cases}
% \end{equation}
%     But
%     \begin{equation}\label{eq:CZ_estimate_thmA}
%         \int_{\B^2}|\nabla u|^p>\Lambda\left(1+ \int_{\B^2}|F|^p\right).
%     \end{equation}
\end{theorem}
%%%%%%%%%%%%%%%%%%%%%%%%%%%%%%%%%%%%%%%%%%%
% \ToDo Examples!

We now explain how \Cref{th:CZ_failure_A}, in particular \Cref{cond:A}, relates to known examples, and provide new ones.
\begin{example}[An example by Johansson -- \cite{Johansson24}]
In this example, along a ray, $A$ asymptotically has a different slope at $+\infty$ than at $-\infty$, with resulting non-uniqueness for very weak solutions. This example falls within the realm of \Cref{th:CZ_failure_A} and thus exhibits failure of the Calder\'on--Zygmund estimate.
 
Let $A:\R^2\to \R^2$ be such that there exists $v_0\in \S^1$ satisfying
$$\lim_{t\to \infty}\frac{A(tv_0)}{t}= a v_0\quad \text{and} \quad \lim_{t\to \infty}\frac{A(-tv_0)}{t}= -b v_0$$
where $0<a<b$. Then, $v_1=v_2=\frac{1}{2}v_0$, $v_3=-v_0$, $w_1=w_2=\frac{a}{2}v_0$, $w_3=-bv_0$ satisfy \Cref{cond:A} with $\beta =a/b$.
\end{example}

\begin{example}
% One of the main motivations of the present work was to generalize the example in \cite{Johansson24,Vincent25} to understand also situations where along any ray the asymptotics of $A$ at $\pm \infty$ are the same, but there is fluctuation in angular direction:
One of the main motivations of the present work was to generalize the example in \cite{Johansson24,Vincent25} to also understand situations where, along any ray, the asymptotics of $A$ at $\pm \infty$ are the same, but there is fluctuation in the angular direction.

    Let $A:\R^2\to \R^2$ be given by 
    $$A(\xi):=\begin{cases}
        \left(2+ \frac{(\xi_1)^2-(\xi_2)^2}{|\xi|^2}\right)\xi, \qquad &\xi=(\xi_1,\xi_2)\neq 0\\
        0,\qquad &\xi=0.
    \end{cases}$$
    Taking 
    $v_1= \frac{1}{2}(1,1),\;v_2= \frac{1}{2}(1,-1),\; v_3=(-1,0),$ $w_1= (1,1),\; w_2=(1,-1),\;w_3=3(-1,0),\; \beta=\frac{2}{3},$
    gives \Cref{cond:A}.
\end{example}

The asymptotic assumptions \eqref{eq:limit of A} and \eqref{eq:three_ray_conv} are stronger than they need to be; for example, they do not recover the example from \cite{Vincent25}. Indeed, it is enough that they hold on sequences with controlled spacing going to infinity.

\begin{example}[Degeneracy along sequences]
For $0<\theta<1$, let $A:\R^2\to \R^2$ be given as 
    \begin{equation*}\label{eq: for A}
        A(\xi):= \begin{pmatrix}
        1-\sqrt{\theta}\sin^2(\pi \xi_1)&0\\
        0&1-\sqrt{\theta}\sin^2(\pi \xi_2)
    \end{pmatrix} \begin{pmatrix}
        \xi_1\\\xi_2
    \end{pmatrix}
    \end{equation*}
    for $\xi=(\xi_1,\xi_2)\in \R^2$. In this example, $A$ satisfies limits as in \Cref{cond:A} along sequences and the spacing between terms grows at a specific rate \cite{Vincent25}.
\end{example}

\subsection*{AI Acknowledgment}
AI use is acknowledged.

\subsection*{Funding Acknowledgment} Both authors received funding through NSF Career DMS-2044898.

\section{Preliminaries}
Let $\mathcal{M}(\R^{2\times 2})$ be the space of (signed) Radon measures on $\R^{2\times 2}$, and let $\mathcal{P}(\R^{2\times 2}) \subset \mathcal{M}(\R^{2\times 2})$ denote the subset of probability measures. We recall the notion of laminates of finite order (see, for example, \cite{AFS08}, \cite{KMSX24}), introduced by Pedregal \cite{Pedregal_1993} and noted earlier by Dacorogna \cite{Dac85}.

\begin{definition}[Elementary splitting and laminates of finite order] Given probability measures $\nu,\mu\in \mathcal{P}(\R^{2\times 2})$ we say that $\mu$ is obtained from $\nu$ by \textit{elementary splitting} if $\nu $ has the form $\nu=\lambda \delta_B +(1-\lambda)\Tilde{\nu}$ for some $\Tilde{\nu}\in \mathcal{P}(\R^{2\times 2})$, there exist matrices $B_1,B_2\in \R^{2\times 2}$ and $\lambda'\in (0,1)$ such that $B=\lambda'B_1+(1-\lambda')B_2$ and rank$(B_2-B_1)=1$, and moreover
$$\mu= \lambda(\lambda'\delta_{B_1}+(1-\lambda')\delta_{B_2})+(1-\lambda)\Tilde{\nu}.$$
The set $\mathcal{L}(\R^{2\times 2})$ of laminates of finite order is defined as the smallest set which is invariant under elementary splitting and contains all  Dirac masses. 

\end{definition}
Let $\nu\in \mathcal{L}(\R^{2\times 2})$ be supported on a finite set of matrices, that is, of the form $\nu= \sum_{j=1}^N \lambda_j \delta_{X_j}$. The center of mass, or \textit{barycenter}, of $\nu$ will be denoted by $\overline{\nu}:= \sum_{j=1}^N \lambda_j X_j$. It is easy to see that the center of mass is invariant under splitting. Let $B, B_1, B_2 \in \R^{2\times 2}$ and $\lambda \in (0,1)$. We say that
$
B = \lambda B_1 + (1-\lambda) B_2
$
is an \textit{elementary rank-one splitting of matrices} if $\operatorname{rank}(B_1 - B_2) = 1$. Laminates of finite order are the basic building blocks for a more refined class of measures, namely \emph{staircase laminates}, which were introduced by Faraco \cite{faracomilton} and have since found various applications in the literature \cite{AFS08, bsv13, colombo2022nonclassicalsolutionsplaplaceequation, CFM05, CFMM05, fls21, FMCO18, KSM03, KMSX24}. 
%For the precise construction of staircase laminates, we refer to \cite[Proposition 3.1]{KMSX24}.

A set $\Omega\subset \R^2$ is called a \textit{regular domain} if it is open, bounded, and connected, and its boundary $\partial\Omega$ has zero $2$-dimensional Lebesgue measure. For $X\in \R^{2\times 2}$ and $b\in \R^2$, we write $l_{X,b}$ for the affine map $l_{X,b}(x)=Xx+b$. The definition below follows the presentation in \cite{KMSX24}.

\begin{definition}[Piecewise affine maps]\label[definition]{defn piecewise affine}
    Let $\Omega$ be a regular domain. We call a map $w\in W^{1,1}(\Omega,\R^2)$ piecewise affine if there exists an at most countable collection of pairwise disjoint regular domains $\Omega_i\subset \Omega$ and a null set $\mathcal{N}$ such that $\Omega=\bigcup_i \Omega_i\cup \mathcal{N}$ and $w$ is an affine map on each $\Omega_i$. That is, for each $i$, there exist $X_i\in \R^{2\times 2}$ and $b_i\in \R^2$ such that $w=l_{X_i,b_i}$ on $\Omega_i$. We also denote by $\dot{\Omega}_w=\bigcup_i\Omega_i$ (or simply $\dot{\Omega}$) the open subset of $\Omega$ where $w$ is locally affine. 
\end{definition}

Next, we state an important proposition about laminates of finite order. There are various versions of the following proposition (see, for example, \cite{Cel93, DM99, DMP08, kir03, KSM03, MS03, Pom10, syc01}). We refer to the proof presented in \cite[Lemma 2.2]{KMSX24}.
%%%%%%%%%%%%%%%%%%%%%%%%%%%%%%%%%%%%%%%%%%%%%%%%
%%%%%%%%%%%%%%%%%%%%%%%%%%%%%%%%%%%%%%%%%%%%%%%%
\begin{proposition}\label{Prop_finite_prelim}
    Let $\nu\in \mathcal{L}(\R^{2\times 2})$ be a laminate of finite order with center of mass $X$. Write $\nu= \sum_{j=1}^J\lambda_j \delta_{X_j}$ with $\lambda_j>0$ and $X_j\neq X_k$ for $j\neq k$. For any $b\in \R^2$, any $\varepsilon>0$ and any regular domain $\Omega\subset \R^2$ there exists a piecewise affine Lipschitz map $w:\Omega\to \R^2$ with
    \begin{itemize}
        \item[(i)] $w= l_{X,b}$ on $\partial \Omega$,
        \item[(ii)] $\Vert \nabla w\Vert_{L^\infty(\Omega)}\leq \max_{j}|X_j|$ 
        \item[(iii)] and 
        \begin{equation}\label{eq:prop prelim}
             (1-\varepsilon)\lambda_j |\Omega| \leq |\{ x\in \Omega:\nabla w=X_j \} |\leq (1+\varepsilon) \lambda_j |\Omega|
        \end{equation}
        for each $j=1,...,J$.
    \end{itemize}
\end{proposition}
%%%%%%%%%%%%%%%%%%%%%%%%%%%%%%%%%%%%%%%%%%
%%%%%%%%%%%%%%%%%%%%%%%%%%%%%%%%%%%%%%%%
\section{Proofs of \Cref{thm:coercive counterexample} and \Cref{th:CZ_failure_A}}
We first show that the condition $H(v)\cdot v \neq 0$ for $v\in \S^1$ 
%in \cref{item 2 them main} of \Cref{th:main} 
is necessary for Calder\'on--Zygmund estimate in \eqref{eq:CZ est main thm} of %\cref{item 1:thm main} of 
\Cref{th:main} to hold. 
%regardless of whether the limit map $H$ in \eqref{eq:limit of A} has a linear extension or not.
%%%%%%%%%%%%%%%%%%%%%%%%%%%%%%%%%%%%%%%%%%%%%%%%%%%%%%%%%%
\begin{proof}[Proof of \Cref{thm:coercive counterexample}]
    %Let $\Tilde{H}$ be the extension of $H$ to $\R^2$ as given in \eqref{eq:extnesion of H}. 
    Suppose there exists a vector $v_0\in \S^1$ such that $H(v_0)\cdot v_0=0$. This implies that there exists a $c_0\in \R$ such that $H(v_0)=c_0 J v_0$ where \begin{equation}\label{eq:J}
             J := \begin{pmatrix} 0 & -1 \\ 1 & 0 \end{pmatrix}.
         \end{equation} Fix any $1<p<\infty$ and any $\Lambda>0$. For $\lambda\in (0,1/2)$ (to be chosen later), the following is an elementary rank-one splitting of matrices:
    $$\begin{pmatrix}
        0\\0
    \end{pmatrix}=\lambda \underbrace{ \begin{pmatrix}
        v_0\\0
    \end{pmatrix}}_{=:G}+
    (1-\lambda)\underbrace{\begin{pmatrix}
        \frac{-\lambda}{1-\lambda} v_0\\0
    \end{pmatrix}}_{=:S}.$$
    For any $t>0$, the following is a laminate of finite order with barycenter $0$:
    $$\lambda \delta_{tG}+(1-\lambda)\delta_{tS}.$$
    For $\varepsilon>0$, by \Cref{Prop_finite_prelim}, there exists a piecewise affine Lipschitz map $$w=(u,\Tilde{u}):\B^2\to \R^2$$ such that $w=0$ on $\partial\B^2$, $\Vert \nabla w\Vert_{L^\infty(\B^2)}\leq t$, and \eqref{eq:prop prelim} holds. On the set $\{x\in \B^2: \nabla w= tG\}$, we have 
    $$|A(\nabla u)-c_0J\nabla u|\leq  |A(tv_0)-tH(v_0)|+|tH(v_0)-c_0J tv_0| < \eta_t t,$$
    for some $\eta_t>0$ and $\lim_{t\to\infty} \eta_t=0$ which comes from \eqref{eq:limit of A}.
    Define $F:= A(\nabla u)-c_0J\nabla u $. Then, 
    $$\div(A(\nabla u))=\div(F)\quad \text{in }\B^2 .$$
    We have 
    \begin{equation}\label{eq:norm F}
        \begin{split}
            \int_{\B^2}|F|^p &= \int_{\{x: \nabla w= tG\}} |F|^p +\int_{\{x: \nabla w= tS\}} |F|^p +\int_{\Omega_{err}}|F|^p\\
            &\aleq  t^p \left( \eta_t^p \lambda+ (1-\lambda)\frac{\lambda^p}{(1-\lambda)^p} +\varepsilon \right)\\
            &\leq C t^p (\eta_t^p +\lambda^p +\varepsilon )
        \end{split}
    \end{equation}
    for a uniform constant $C>0$ (depending only on $p$ and $L$). Now, 
    \begin{equation}\label{eq:norm u}
        \int_{\B^2} |\nabla u|^p \geq \int_{\{x: \nabla w= tG\}} |\nabla u|^p\;\ageq t^p  \lambda \geq c t^p \lambda
    \end{equation}
    where $c>0$ depends only on $p$. 

    Since $p>1$, we first choose $\lambda$ small enough so that 
    $$C\lambda^{p-1} <\frac{c}{6\Lambda}.$$
    Then, choose $t$ large enough so that 
    $$C \eta_t^p <\frac{c\lambda }{6\Lambda}\quad \text{and}\quad ct^p\lambda >2\Lambda. $$
    Finally, we choose $\eps$ small enough such that
\[
C\eps <\frac{c \lambda}{6\Lambda}.
\]
    Combining these with \eqref{eq:norm F} and \eqref{eq:norm u} gives $$\Lambda \left(1+\int_{\B^2} |F|^p \right) < \Lambda + \frac{c\lambda t^p}{2} < \frac{c\lambda t^p}{2}+\frac{c\lambda t^p}{2} \leq \int_{\B^2}|\nabla u|^p$$
    proving \eqref{eq:CZ_estimate_thmA}. This proves the theorem.
\end{proof}
%===============================================================
Next, we prove \Cref{th:CZ_failure_A}. Let $v_i,w_i$ for $i=1,2,3$ be as in \Cref{cond:A}. We define 
\begin{equation}\label{eq:S0 Gi}
    \begin{split}
        S_0&:= \begin{pmatrix}
    -v_3 \\ -w_3^\perp
\end{pmatrix},\\
G_i&:= \begin{pmatrix}
    v_i \\ -w_i^\perp
\end{pmatrix},\quad i=1,2,3,
    \end{split}
\end{equation}
where $v^\perp = J v$ for $v\in \R^2$ with $J$ as in \eqref{eq:J}.
% Let 
% \begin{equation}\label{eq:man/min-norm of vectors}
%     c_1:= \max\{|v_i|,|w_i|;\; i=1,2,3\}\quad\text{and}\quad c_2:= \min\{|v_i|,|w_i|;\; i=1,2,3\}>0.
% \end{equation}
\begin{proposition}\label{prop_basic_splitting}
    Let $S_0$ and $G_i$ for $i=1,2,3$ be as in \eqref{eq:S0 Gi}. There exist numbers $\lambda_i\in (0,1)$ for $i=1,2,3$, $t_1,t_2\in (0,\infty)$, $t_3>1$ and $T_1,T_2\in \R^{2\times 2}$ such that the following are elementary rank-one splittings of matrices:
    \[
    \begin{split}
        S_0&= \lambda_1 t_1 G_1 +(1-\lambda_1)T_1,\\
        T_1&= \lambda_2 t_2 G_2 +(1-\lambda_2)T_2,\\
        T_2&= \lambda_3 t_3 G_3 +(1-\lambda_3)t_3 S_0.\\
    \end{split}
    \]
 Moreover, for $\gamma:=(1-\lambda_1)(1-\lambda_2)(1-\lambda_3)$, we can ensure that
 \begin{equation}\label{eq:bound-on-gamma-t3}
     t_3\gamma <1.
 \end{equation}
 \end{proposition}
 The proof of the above proposition immediately follows from the following lemmas.
%============================================
%Lemma for tau_1 and tau_2
\begin{lemma}\label{lemma_tau_i}
    For $i=1,2$, there exists a unique $\tau_i\in (0,\infty)$ such that $\det (S_0-\tau_i G_i)=0$. 
\end{lemma}
\begin{proof}
    Fix $i=1,2$. For $t\in \R$, we define $$p(t):=\det(S_0-tG_i)=(-v_3-tv_i)\cdot (-w_3+tw_i).$$
    The quadratic coefficient of $p(t)$ is $-v_i\cdot w_i <0$ from \eqref{eq:coercivity}. So, $\lim_{t\to\infty}p(t)=-\infty$. At $t=0$, $p(0)=v_3\cdot w_3>0$. Since $p(t)$ is a quadratic polynomial, $p(t)$ has exactly one positive root and it is a simple root. We call it $\tau_i$.  
\end{proof}
%===============================================
%lemma for t_2
\begin{lemma}\label{lemma_t_2}
    Let $\tau_1$ be as in \Cref{lemma_tau_i}. For any $\lambda\in (0,1)$, define 
    \begin{equation}\label{eq:T_1}
       T_1(\lambda):= \frac{S_0-\lambda \tau_1 G_1}{(1-\lambda)}. 
    \end{equation}
    Then, there exist $\varepsilon\in (0,1)$ and a smooth map $t_2: (0,\varepsilon)\to (0,\infty)$ such that 
    $$\det (T_1(\lambda)-t_2(\lambda) G_2)=0$$
    and $\lim_{\lambda\to 0^+} t_2(\lambda)=\tau_2$.
\end{lemma}
\begin{proof}
    Define 
    $$g(\lambda,t):= \det(T_1(\lambda)- t G_2)$$ for $(\lambda,t)\in (-1,1)\times (0,\infty)$. Clearly, $g$ is a smooth map and $g(0, \tau_2)=0$ for $\tau_2$ as in \Cref{lemma_tau_i}. Moreover, we get $\frac{\partial}{\partial t}g(0, \tau_2)\neq0$ since $\tau_2$ is a simple root of $p(t)$ as in \Cref{lemma_tau_i}. By the implicit function theorem, there exist $0<\varepsilon<1,\;0<\chi<\tau_2 $ and a smooth map 
    $t_2:(0,\varepsilon)\to (\tau_2-\chi , \tau_2 +\chi)$ such that $$g(\lambda, t_2(\lambda))=0$$
    whenever $0<\lambda<\varepsilon$ and clearly, $\lim_{\lambda\to 0^+} t_2(\lambda)=\tau_2$.
\end{proof}
%===================================================
\begin{lemma}\label{lemma:first and second split}
    Let $\tau_i$ be as in \Cref{lemma_tau_i}, $T_1$ be as in \eqref{eq:T_1}, and $t_2:(0,\varepsilon)\to (0,\infty)$ be as in \Cref{lemma_t_2}. Then, there exists $0<\varepsilon_1\ll 1$ such that for any $0<\lambda<\varepsilon_1$, the following are elementary rank-one splittings of matrices:
    \begin{equation}\label{eq:first split}
        S_0= \lambda \tau_1 G_1+(1-\lambda) T_1 (\lambda),
    \end{equation}
    \begin{equation}\label{eq:second split}
        T_1(\lambda)= \lambda_2(\lambda) t_2(\lambda) G_2 + (1-\lambda_2(\lambda)) T_2(\lambda)
    \end{equation}
    where 
    \begin{equation}\label{eq:lambda_2}
        \lambda_2(\lambda):= \frac{\lambda \tau_1}{(1-\lambda)t_2(\lambda)}
    \end{equation}
    and
    \begin{equation}\label{eq:T_2}
        T_2(\lambda):= \frac{T_1(\lambda)-\lambda_2(\lambda)t_2(\lambda)G_2}
{1-\lambda_2(\lambda)}.
    \end{equation}
\end{lemma}
\begin{proof}
    Since $\det(S_0-\tau_1 G_1)=0$ from \Cref{lemma_tau_i} and $S_0-\tau_1 G_1 \neq 0$ from \eqref{eq:coercivity}, we get rank$(S_0-\tau_1 G_1)=1$. So, \eqref{eq:first split} is an elementary rank-one splitting for any $0<\lambda<1$. 
    %Observe that
%$$0<\frac{\lambda \tau_1}{(1-\lambda)t_2(\lambda)}< \frac{\lambda \tau_1}{(1-\lambda)(\tau_2 -\delta)}
%so choose \varepsilon_1= min \varepsilon, \frac{\tau_2-\delta}{\tau_1+\tau_2-\delta}$$
Since $t_2(\lambda)>0$, we have $\lambda_2(\lambda)>0$ for all $0<\lambda<\varepsilon$. Since $t_2(\lambda)$ is smooth and $\lim_{\lambda \to 0^+}t_2(\lambda)=\tau_2>0$, we have $$\lim_{\lambda\to 0^+}\lambda_2(\lambda)=0.$$
Hence, there exists a positive number $\varepsilon_1< \varepsilon $ such that $0<\lambda_2(\lambda)<1$ for all $0<\lambda<\varepsilon_1$. From \Cref{lemma_t_2}, we have $\det(T_1(\lambda)-t_2(\lambda)G_2)=0$ and $T_1(\lambda)-t_2(\lambda)G_2\neq 0$ for sufficiently small $\lambda$. Thus, \eqref{eq:second split} is an elementary rank-one splitting for any $0<\lambda<\varepsilon_1$.
\end{proof}
%====================================================
\begin{lemma}\label{lemma_third split}
    Let $T_2(\lambda)$ be as in \eqref{eq:T_2}. There exists $0<\varepsilon_2<1$ such that for any $0<\lambda<\varepsilon_2$
    \begin{equation}\label{eq:third split}
        T_2(\lambda)
=\lambda_3(\lambda)t_3(\lambda)G_3
+(1-\lambda_3(\lambda))t_3(\lambda)S_0
    \end{equation}
is an elementary rank-one splitting for some $\lambda_3(\lambda)\in (0,1)$ and $t_3(\lambda)>1$. Moreover, we have
$$t_3(\lambda)(1-\lambda)(1-\lambda_2(\lambda))
(1-\lambda_3(\lambda))<1.$$
% More precisely,
% \[
% t_3(\lambda)
% =\frac{1+\beta\lambda t_1}
% {(1-\lambda)(1-\lambda_2(\lambda))}
% \]
% and
% \[
% \lambda_3(\lambda)
% =\frac{(1+\beta)\lambda t_1}
% {2(1+\beta\lambda t_1)}.
% \]
\end{lemma}
\begin{proof}
    From \eqref{eq:lambda_2} and \eqref{eq:T_1}, we get
\[
\begin{split}
    T_2(\lambda)
&=\frac{S_0-\lambda \tau_1(G_1+G_2)}
{(1-\lambda)(1-\lambda_2(\lambda))}\\
&= \frac{1}{(1-\lambda)(1-\lambda_2(\lambda))}
\begin{pmatrix}
-v_3-\lambda \tau_1(v_1+v_2)\\
-w_3^\perp+\lambda \tau_1(w_1^\perp+w_2^\perp)
\end{pmatrix}.
\end{split}
\]
Using \eqref{eq:unit_vector_sum_v} and \eqref{eq:unit_vector_sum_w},
% \[
% v_1+v_2=-v_3,
% \qquad
% w_1^\perp+w_2^\perp=-\beta w_3^\perp,
% \]
we have
\begin{equation}\label{eq:forLHS 3trd splitting}
   \begin{aligned}
T_2(\lambda)
=\frac{1}{(1-\lambda)(1-\lambda_2(\lambda))}
\begin{pmatrix}
(\lambda \tau_1-1)v_3\\
-(1+\beta\lambda \tau_1)w_3^\perp
\end{pmatrix}.
\end{aligned} 
\end{equation}
We want to find $\lambda_3\in(0,1)$ and $t_3>1$
such that
\[
T_2(\lambda)
=\lambda_3t_3G_3+(1-\lambda_3)t_3S_0.
\]
Observe that
\begin{equation}\label{eq:RHS for third splitting}
    \begin{aligned}
\lambda_3t_3G_3+(1-\lambda_3)t_3S_0
&=\lambda_3t_3
\begin{pmatrix}
v_3\\
-w_3^\perp
\end{pmatrix}
+(1-\lambda_3)t_3
\begin{pmatrix}
-v_3\\
-w_3^\perp
\end{pmatrix}\\
&=
\begin{pmatrix}
(2\lambda_3-1)t_3v_3\\
-t_3w_3^\perp
\end{pmatrix}.
\end{aligned}
\end{equation}
Comparing the second rows of \eqref{eq:forLHS 3trd splitting} and \eqref{eq:RHS for third splitting}, we choose
\begin{equation}\label{eq:t_3}
    t_3(\lambda)
=\frac{1+\beta\lambda \tau_1}
{(1-\lambda)(1-\lambda_2(\lambda))}.
\end{equation}
Comparing the first rows of \eqref{eq:forLHS 3trd splitting} and \eqref{eq:RHS for third splitting} gives
\[
(2\lambda_3-1)t_3
=\frac{\lambda \tau_1-1}
{(1-\lambda)(1-\lambda_2(\lambda))}.
\]
Substituting the expression for $t_3$, we get
\[
2\lambda_3-1
=\frac{\lambda \tau_1-1}{1+\beta\lambda \tau_1}.
\]
Thus,
\begin{equation}\label{eq:lambda_3}
    \lambda_3(\lambda)
=\frac{(1+\beta)\lambda \tau_1}
{2(1+\beta\lambda \tau_1)}.
\end{equation}

It remains to verify that
$0<\lambda_3(\lambda)<1$ and $t_3(\lambda)>1$.

Since 
$0<\lambda, \lambda_2(\lambda)<1$, we have 
\[
t_3(\lambda)
=\frac{1+\beta\lambda \tau_1}
{(1-\lambda)(1-\lambda_2(\lambda))}>1.
\]
Moreover, since $\beta>0$ and $\tau_1>0$,
\[
\lambda_3(\lambda)
=\frac{(1+\beta)\lambda \tau_1}
{2(1+\beta\lambda \tau_1)}>0.
\]
Since $\lambda_3(\lambda)>0$ and $\lim_{\lambda\to 0^+} \lambda_3(\lambda)=0$,
there exists $\varepsilon_2>0$ such that
\[
0<\lambda_3(\lambda)<1\quad\forall\;0<\lambda<\varepsilon_2.
\]

Finally,
\[
t_3(\lambda)G_3-t_3(\lambda)S_0
=t_3(\lambda)
\begin{pmatrix}
2v_3\\
0
\end{pmatrix},
\]
which has rank one. Therefore, \eqref{eq:third split} is an elementary rank-one splitting.

Finally, from \eqref{eq:lambda_3} and \eqref{eq:t_3}, we have
\[
\begin{aligned}
t_3(\lambda)(1-\lambda)(1-\lambda_2(\lambda))
(1-\lambda_3(\lambda))
&=(1+\beta\lambda \tau_1)(1-\lambda_3(\lambda))\\
&=(1+\beta\lambda \tau_1)
\left(1-\frac{(1+\beta)\lambda \tau_1}
{2(1+\beta\lambda \tau_1)}\right)\\
% &=1+\beta\lambda t_1
% -\frac{(1+\beta)\lambda t_1}{2}\\
&=1+\frac{\beta-1}{2}\lambda \tau_1<1,
\end{aligned}
\]
where the last inequality follows from
$0<\beta<1$, $\tau_1>0$, and $\lambda>0$.
\end{proof}
%================================================
\begin{proof}[Proof of \Cref{prop_basic_splitting}]
    Choose $\lambda>0$ sufficiently small so that we can apply \Cref{lemma_t_2}, \Cref{lemma:first and second split} and \Cref{lemma_third split} simultaneously. Take $\lambda_1=\lambda$ and $t_1=\tau_1$. This gives the desired elementary splittings with the coefficients.
\end{proof}
%================================================

Now, we define the threshold exponent as follows:
\begin{equation}\label{eq:thershold}
    \bar{q}:=\frac{-\log \gamma}{\log t_3}
\end{equation}
where $\gamma \in (0,1)$ and $t_3>1$ are as in \Cref{prop_basic_splitting}. From \eqref{eq:bound-on-gamma-t3}, we get $$
    \log (t_3 \gamma)<0
    \iff  1<\frac{-\log \gamma }{\log t_3}.$$
Thus, $\bar{q}>1$.

% For $t_i$ as in \Cref{prop_basic_splitting}, let
% \begin{equation}\label{eq:man/min-norm of vectors}
%     c_1:= \max\{|t_iv_i|,|t_i w_i|;\; i=1,2,3\}\quad\text{and}\quad c_2:= \min\{|t_iv_i|,|t_i w_i|;\; i=1,2,3\}>0.
% \end{equation}
%=====================================
% \begin{lemma}\label{lemma:q_bar}
%     Let $\bar{q}$ be as in \eqref{eq:thershold}. Then, $\bar{q}>1$. Also, for any $1<p<\bar{q}$, we have 
% \end{lemma}
%=================================================
\begin{proposition}\label{prop_laminate_construction}
    Let $S_0$ and $G_i$ for $i=1,2,3$ be as in \eqref{eq:S0 Gi}. Let $\lambda_i, \;t_i$ for $i=1,2,3$, and $\gamma$ be as in \Cref{prop_basic_splitting}. For every $N \geq 1$,
    \[
\begin{split}
\nu_N := \frac{1}{2} \delta_{G_3}+&\frac{1}{2}\sum_{k=0}^{N-1} \gamma^k \Big[\, \lambda_1 \,\delta_{t_3^k t_1 G_1}
  + (1-\lambda_1)\lambda_2 \,\delta_{t_3^k t_2 G_2} 
  + (1-\lambda_1)(1-\lambda_2)\lambda_3 \,\delta_{t_3^{k+1} G_3} \Big]
  \\+ & \frac{1}{2}\gamma^N \delta_{t_3^N S_0}
\end{split}
\]
is a laminate of finite order with barycenter
\[
\begin{pmatrix} 0 \\ -w_3^{\perp} \end{pmatrix}.
\]
% is a laminate of finite order with barycenter $S_0$. 
Moreover, for any $1<p<\bar{q}$ for $\bar{q}$ as in \eqref{eq:thershold}, we have  
 \begin{equation}\label{eq:bound-on-gamma-t3-p}
     t_3^p\gamma <1.
 \end{equation}
%  In addition, 
%  \[
% \Tilde{\nu}_N := \frac{1}{2}\,\delta_{G_3} + \frac{1}{2}\,\nu_N
% \]
% is a laminate of finite order with barycenter
% \[
% \begin{pmatrix} 0 \\ -w_3^{\perp} \end{pmatrix}.
% \]
 \end{proposition}
\begin{proof}
    Starting with $S_0$, we apply the three elementary rank-one splittings in \Cref{prop_basic_splitting}. At the next step, we apply the same splittings, multiplied with $t_3$, to the matrix $t_3 S_0$.  For any $N\geq 1$, we repeat this $N$ times to get a laminate of finite order with barycenter $S_0$. For the initial splitting, we have 
$$\begin{pmatrix} 0 \\ -w_3^{\perp} \end{pmatrix}=\frac{1}{2}\underbrace{\begin{pmatrix} v_3 \\ -w_3^{\perp} \end{pmatrix}}_{G_3}+\frac{1}{2}\underbrace{\begin{pmatrix} -v_3 \\ -w_3^{\perp} \end{pmatrix}}_{S_0}.$$ 
Combining these gives the laminate $\nu_N$. Next, let $1<p<\bar{q}$.
%\akshara{Looks like this argument works for any $0<p<\bar{q}$}. 
Since $t_3>1$,  we have 
\[
\begin{split}
    p<\frac{-\log \gamma}{\log t_3}
    % \iff \log t_3^p &<-\log \gamma\\
    \iff \log (t_3^p \gamma)<0
    \iff  \;t_3^p\gamma <1.
\end{split}
\]
This completes the proof.
\end{proof}
%=====================================
%%%%%%%%%%%%%%%%%%%%%%%%%%%%%%%%%%%%%%%%%%
\begin{proof}[Proof of \Cref{th:CZ_failure_A}]
    Let $\bar{q}$ be as in \eqref{eq:thershold}. Fix $\Lambda>0$ and $1<p<\bar{q}$. Let $N\geq 1$ and $R>0$ (to be chosen depending on $\Lambda$). We apply \Cref{prop_laminate_construction} to find a laminate of finite order $\nu_N$. Multiply every matrix in the support of $\nu_N$ by $R$ to obtain the laminate $\Tilde{\nu}_N$. Let $\varepsilon\in (0,1/2)$ (to be chosen small enough depending on $N$ and $R$). Next, we apply \Cref{Prop_finite_prelim} to $\Tilde{\nu}_N$ to obtain a piecewise affine Lipschitz map 
    $$w=(u,\Tilde{u}):\B^2\to \R^2$$ such that 
    \begin{equation}\label{eq:boundary_data_thmA}
        w(x)= R\begin{pmatrix}
        0\\ -w_3^\perp 
    \end{pmatrix}x\quad \text{on}\;\partial\B^2
    \end{equation}
    and $$\Vert \nabla w\Vert_{L^\infty(\B^2)}\aleq c_1 R t_3^N$$
    for $c_1:= \max\{t_i^\ell|v_i|^m,t_i^\ell |w_i|^m:\; i=1,2,3;\;\ell,m=0,1\}>0$. 
    
    For $S\in \R^{2\times 2}$, we denote $\Omega_S:= \{x\in \B^2:\;\nabla w(x)=S\}$. Let 
    $$\Sigma_i:=\{t_3^k t_i:\ k=0,\dots,N-1\},\;\;i=1,2,\quad \Sigma_3:= \{t_3^k:\ k=0,\dots,N\},$$
    $$\mathcal{G}:= \{RsG_i:\;i=1,2,3, \; s\in \Sigma_i\}$$
    and 
    $$\Omega_{err}:=\B^2 \backslash \left(\Omega_{Rt_3^N S_0} \dot{\cup} \dot{\bigcup_{S\in \mathcal{G}}}\Omega_S\right)$$
    where $\dot{\cup}$ denotes a disjoint union of sets. From \Cref{Prop_finite_prelim}, 
    we get $$|\Omega_{err}|\leq \varepsilon |\B^2|.$$ Set $F:=A(\nabla u)-J\nabla \Tilde{u}.$
    \begin{equation}\label{eq:defn_of_FthmA}
       \div(A(\nabla u))=\div(F)\quad \text{in}\;\B^2.
    \end{equation}
     We have
    \begin{equation}\label{eq: f bound on err set}
        \int_{\Omega_{err}}|F|^p\aleq_A \varepsilon \Vert \nabla w\Vert^p_{L^{\infty}(\B^2)} |\B^2|\leq C_1\varepsilon R^p t_3^{Np}
    \end{equation}
    where $C_1>0$ depends only on $A$ and $p$. 
    
     Next, we define 
    $$\eta_{R,N}:=\max \left\{ \left|\frac{A(Rs v_i)}{Rs}-w_i\right|:\quad i=1,2,3,\;s\in \Sigma_i  \right\}.$$
    From \eqref{eq:three_ray_conv}, we get
    \begin{equation}\label{eq:convergence of eta}
        \eta_{R,N}\to 0\quad \text{as}\quad R\to \infty
    \end{equation}
    for each fixed $N$.
%     Next, take any $G\in \mathcal{G}$, say 
%     $$S= Rt_3^k t_i \begin{pmatrix}
%         v_i\\-w_i^\perp 
%     \end{pmatrix}.$$
% Then, for any $x\in\Omega_G$, we have
% $$|F(x)|=|A(Rt_3^k t_i v_i)-Rt_3^k t_i w_i|\leq c_1 \eta_{R,N} R t_3^k.$$
Take any $S\in\mathcal G$, say $S=Rs\,G_i$ with $s\in\Sigma_i$. Then $\nabla u=Rsv_i$ and
$J\nabla\tilde u=Rsw_i$ on $\Omega_S$, so
\[
|F|=|A(Rsv_i)-Rsw_i|\le \eta_{R,N}\,Rs\le c_1\,\eta_{R,N}\,R\,t_3^{k}
\quad\text{on }\Omega_S,
\]
where $k\in\{0,\dots,N\}$ is the power of $t_3$ in $s$.
 From \eqref{eq:prop prelim} and \eqref{eq:bound-on-gamma-t3-p}, we get
\begin{equation}\label{eq: f bound on G}
    \sum_{G\in \mathcal{G}} \int_{\Omega_G}|F|^p \aleq c_1^p  |\B^2| \eta_{R,N}^p R^p\sum_{k=0}^{N}(\gamma t_3^p)^k
\le C_2 \eta_{R,N}^pR^p
\end{equation}
where $C_2>0$ depends only on $A$ and $p$. Next, applying \eqref{eq:prop prelim} on $\Omega_{Rt_3^N S_0}$, we get
\begin{equation}\label{eq: f bound on S_N}
    \int_{\Omega_{R t_3^N S_0}}|F|^p \leq C_3 R^p (\gamma t_3^p)^N 
\end{equation}
where $C_3>0$ depends only on $A$ and $p$. Therefore, \eqref{eq: f bound on err set}, \eqref{eq: f bound on G} and \eqref{eq: f bound on S_N} together give
\begin{equation}\label{eq:bound on F}
\int_{\B^2}|F|^p
\le CR^p\left[\eta_{R,N}^p+(\gamma t_3^p)^N+\eps t_3^{Np}\right]
\end{equation}
for some constant $C>0$ depending only on $A$ and $p$. We observe that
\begin{equation}\label{eq:bound lower of nabla u}
    \int_{\B^2}|\nabla u|^p \geq \int_{\Omega_{RG_3}}|\nabla u|^p \geq (R |v_3|)^p |\Omega_{RG_3}|\geq c R^p 
\end{equation}
for some $c>0$ depending only on $A$ and $p$. 

Now, we choose $N$ large enough such that
\[
C(\gamma t_3^p)^N<\frac{c}{6\Lambda}.
\]
Using \eqref{eq:convergence of eta}, we can choose $R$ large enough such that
\[
C\eta_{R,N}^p<\frac{c}{6\Lambda}
\qquad \text{and} \qquad
cR^p>2\Lambda.
\]
Finally, choose $\eps$ small enough such that
\[
C\eps t_3^{Np}<\frac{c}{6\Lambda}.
\]
Combining \eqref{eq:bound on F} and \eqref{eq:bound lower of nabla u} gives
$$\Lambda \left(1+\int_{\B^2} |F|^p \right) < \Lambda + \frac{cR^p}{2} < \frac{cR^p}{2}+\frac{cR^p}{2} \leq \int_{\B^2}|\nabla u|^p.$$
This gives \eqref{eq:CZ_estimate_thmA}. From \eqref{eq:boundary_data_thmA} and \eqref{eq:defn_of_FthmA}, we get $u$ and $F$ satisfy \eqref{eq:PDE_thmA}. This completes the proof.
\end{proof}
%%%%%%%%%%%%%%%%%%%%%%%%%%%%%%%%%%%%%%%%%%
\section{Proof of \Cref{th:main}}
%=====================================
The following proposition gives an equivalence of the linearity of the extension map $\Tilde{H}$ in \eqref{eq:extnesion of H} with \Cref{cond:A}, provided it is non-degenerate. 
\begin{proposition}\label{prop 3vectors}
    Let $B:\R^2\to \R^2$ be continuous, positively homogeneous and $B(v)\cdot v\neq 0$ for $v\neq 0$. Then the following statements are equivalent:
    \begin{enumerate}
        \item \label{item:first} $B$ is linear.
        \item \label{item:scond} Whenever $v_1\neq 0$, $v_2\neq 0$ and $v_1+v_2\neq 0$, we have 
        \begin{equation}\label{eq:Theta_H}
            \Theta(v_1,v_2):= \frac{\big(B(v_1)+B(v_2)\big)\cdot (v_1+v_2)}{B(-v_1-v_2)\cdot (-v_1-v_2)}\notin (0,1).
        \end{equation}
         \item \label{item:thrid} For every $c\in \R$, there are no non-zero vectors $v_1,v_2,v_3\in \R^2$ and $\theta\in (0,1)$ such that 
         \begin{equation}\label{eq:cond A for A+cJ}
             v_1+v_2+v_3=0\quad \text{and}\quad (B+cJ)(v_1)+(B+cJ)(v_2)+\theta (B+cJ)(v_3)=0
         \end{equation}
         for $J$ as in \eqref{eq:J}.
    \end{enumerate} 
    \end{proposition}
    \begin{proof}
        First, we prove \cref{item:first} and \cref{item:scond} are equivalent. Suppose $B$ is linear. Then, $\Theta(v_1,v_2)=1$ whenever $v_1\neq 0$, $v_2\neq 0$ and $v_1+v_2\neq 0$. This implies \cref{item:scond}. 
        
        For the other direction, assume \cref{item:scond}. Define $f:\R^2\to \R$ by 
        \begin{equation}\label{eq:H(x)dot x}
            f(v):= B(v)\cdot v.
        \end{equation}
        Clearly, $f$ is continuous and $f(x)\neq 0$ on $\R^2\backslash\{0\}$. Since $\R^2\backslash\{0\}$ is path connected and $f$ is continuous, $f(\R^2\backslash \{0\})$ is also path connected. Therefore, $f$ has constant sign on $\R^2\backslash \{0\}$. Without loss of generality, assume $f(v)>0$ whenever $v\neq 0$. (Otherwise, we do the same argument for $-B$ and $\Theta$ remains the same for both.)

        From \eqref{eq:Theta_H} and positive homogeneity of $B$, we get
        $$\Theta(v,v)=\frac{f(v)}{f(-v)}\geq 1 \quad \text{and}\quad \Theta(-v,-v)=\frac{f(-v)}{f(v)}\geq 1.$$
        This implies $f(v)=f(-v)$ for $v\neq 0$ and hence, $f$ is an even function.

        Fix $v_1\neq 0$. The function $v_2 \mapsto \Theta(v_1,v_2)$ is a continuous map on the connected set $\R^2\backslash\{0,-v_1\}$ and tends to $1$ as $v_2\to 0$. From \eqref{eq:Theta_H}, $\Theta(v_1,\cdot )$ takes values in $[1,\infty)$ in a small neighborhood of $v_2=0$. By the continuity of the map and connectedness of the domain, we get $\Theta(v_1,v_2)\geq 1$ for all $v_2\in \R^2\backslash\{0,-v_1\}$. This happens for all $v_1\neq 0$. In particular, using the fact $f$ is even, we get
        \begin{equation}\label{eq:inequality of H}
            \left( B(v_1)+B(v_2)\right) \cdot (v_1+v_2) \geq B(v_1+v_2)\cdot (v_1+v_2) \quad \forall\;v_1,v_2\in \R^2.
        \end{equation}

        First, we prove that $B$ is an odd function. Indeed, fix $v_1\neq 0$ and take any $v_2\neq 0$. For $t>0$, using \eqref{eq:inequality of H} we have
        \begin{equation}\label{eq:proving H odd}
            \left( B(v_1)+B(-v_1+tv_2)\right) \cdot t v_2 \geq B(tv_2)\cdot tv_2= t^2 B(v_2)\cdot v_2.
        \end{equation}
        Canceling $t$ from both sides and letting $t\to 0^+$, we get
        $$\left( B(v_1)+B(-v_1)\right)\cdot v_2\geq 0.$$
        Repeating the same argument with $-v_2$ gives 
        $$\left( B(v_1)+B(-v_1)\right)\cdot v_2\leq 0.$$ 
        This happens for any $v_2\neq 0$ and hence, we get $ B(v_1)=-B(-v_1)$ for $v_1\neq 0$. Thus, $ B$ is odd. In addition, using the positive homogeneity and continuity of $B$, we get 
        \begin{equation}\label{eq:mutiplicativity_of H}
            B(tv)=tB(v)\quad\forall\;v\in \R^2,\;t\in \R.
        \end{equation}

        To prove the additivity of $B$, fix any linearly independent vectors $v_1,v_2\in \R^2$. Let $\Delta:= B(v_1)+B(v_2)-B(v_1+v_2)$. From \eqref{eq:inequality of H}, we get
        $$\Delta\cdot (v_1+v_2)\geq 0.$$
        Applying \eqref{eq:inequality of H} for the vectors $v_1+v_2$, $-v_2$ and the vectors $v_1+v_2$, $-v_1$ we get 
        \begin{equation}\label{eq:Delta x and y}
            \Delta\cdot v_1 \leq0 \quad \text{and} \quad \Delta\cdot v_2\leq 0,
        \end{equation}
        respectively. Thus, $\Delta \cdot (v_1+v_2)= 0$. From \eqref{eq:Delta x and y}, $\Delta \cdot v_1=0$ and $\Delta\cdot v_2 =0$. Since $v_1$ and $v_2$ are linearly independent, this implies $\Delta=0$. This implies 
        \begin{equation}\label{eq:additivity of H}
            B(v_1+v_2)=B(v_1)+B(v_2)
        \end{equation}
        for any linearly independent vectors $v_1,v_2\in \R^2$.
Combining \eqref{eq:mutiplicativity_of H} and \eqref{eq:additivity of H}, we get $B$ is linear.

        Next, we prove \cref{item:scond} and \cref{item:thrid} are equivalent. Suppose \cref{item:scond} is not true. Then there exist non-zero vectors $v_1, v_2\in \R^2$ such that $v_1+v_2\neq 0$ and $\Theta(v_1,v_2)\in (0,1)$. Let $v_3:=-(v_1+v_2)$ and $\theta:= \Theta(v_1,v_2)$. Then,
        $$\left(B(v_1)+B(v_2) + \theta B(v_3)\right)\cdot v_3 =0.$$
        Thus, there exists some $\mu\in \R$ such that $$B(v_1)+B(v_2) + \theta B(v_3)= \mu J v_3.$$ Let $c:= \mu/ (1-\theta)$. Then 
        $$B(v_1)+B(v_2) + \theta B(v_3)+c(\theta-1) Jv_3 =0.$$
        This gives 
        $$(B+cJ)(v_1)+(B+cJ)(v_2)+\theta (B+cJ)(v_3)=0.$$
        This implies \cref{item:thrid} fails.
        
        For the other direction, suppose \cref{item:thrid} is not true. Then there exist $c\in \R$, nonzero vectors $v_1,v_2,v_3\in \R^2$, and $\theta \in (0,1)$ such that $v_1+v_2+v_3=0$ and 
        $$(B+cJ)(v_1)+(B+cJ)(v_2)+\theta (B+cJ)(v_3)=0.$$ Taking the scalar product with $v_3=-(v_1+v_2)$ in the above equation gives 
        $$\Theta(v_1,v_2)=\theta \in (0,1).$$ This implies \cref{item:scond} fails. This completes the proof.
    \end{proof}
%%%%%%%%%%%%%%%%%%%%%%%%%%%%%%%%%%%%%%%%%%%%%%%%%%%%%
\begin{proof}[Proof of \Cref{th:main}]
    Let $\Tilde{H}$ be as defined in \eqref{eq:extnesion of H}. Clearly, this is a well-defined extension of the map $H$ to $\R^2$. From \eqref{eq:limit of A}, we get $H$ is continuous (since it is the uniform limit of a continuous map), and so is its positively homogeneous extension $\Tilde{H}$ (follows from \eqref{eq:extnesion of H}). In addition, using \eqref{eq:linearity of A} and \eqref{eq:limit of A}, we have
    \begin{equation}\label{eq:linearity of H}
        |\Tilde{H}(v)|\leq  L|v|\quad \forall\; v\in \R^2.
    \end{equation}
    
    First, we prove \cref{item 2 them main} implies \cref{item 1:thm main}. Assume $H(v)\cdot v\neq 0$ for $v\in \S^1$ and $\Tilde{H}$ is linear. This implies $\Tilde{H}(v)\cdot v\neq 0$ for $v\neq 0$. Without loss of generality, we may assume $\Tilde{H}(v)\cdot v>0$ on $\R^2\backslash \{0\}$ by a similar argument as in \eqref{eq:H(x)dot x}. Since $v\mapsto \Tilde{H}(v)\cdot v$ is a continuous map, there exists $\iota>0$ such that $\Tilde{H}(v)\cdot v\geq \iota$ on $\S^1$. By linearity of $\Tilde{H}$, we get 
    \begin{equation}\label{eq:coercivity of H}
        \Tilde{H}(v)\cdot v\geq \iota |v|^2\quad \forall \;v\in \R^2.
    \end{equation}
    From \eqref{eq:limit of A} and linearity of $\Tilde{H}$, for $\varepsilon\in (0,1)$, there exists $t_\varepsilon >0$ such that 
    $$|A(v)- \Tilde{H}v|<\varepsilon |v| \quad \forall\;|v|>t_\varepsilon.$$
    Since $\Tilde{H}$ is continuous, 
    $c_{\varepsilon}:= \max\limits_{|v|\leq t_\varepsilon}|A(v)-\Tilde{H}v|<\infty$. Thus, 
    \begin{equation}\label{eq: bonund on A-H}
        |A(v)- \Tilde{H}v|<\varepsilon |v| +c_\varepsilon \quad \forall\;v\in \R^2.
    \end{equation}
    Let $1<p<\infty$ and let $u\in W^{1,p}_0(\B^2)$ and $F\in L^p(\B^2,\R^2)$ be such that
    $$\div(A(\nabla u))=\div (F)\quad \text{in}\;\B^2.$$ We look at the equation 
    $$\div(\Tilde{H}\nabla u)=\div(F+\Tilde{H}\nabla u -A(\nabla u)).$$
    Using \eqref{eq:linearity of H} and \eqref{eq:coercivity of H}, and the standard $L^p$-estimates for elliptic equations (for instance, see \cite[Chapter 7]{GiaquintaMartinazzi2012}), we get 
    $$\Vert \nabla u\Vert_{L^p(\B^2)}\leq C_{\Tilde{H},p}\left(\Vert F\Vert_{L^p(\B^2)}+\varepsilon \Vert \nabla u\Vert_{L^p(\B^2)}+c_\varepsilon  |\B^2|^{1/p}\right).$$
    We choose $\varepsilon>0$ small enough so that $C_{\Tilde{H},p} \varepsilon< 1/2$ and absorb the term $\varepsilon \Vert \nabla u\Vert_{L^p(\B^2)}$ to the left-hand side since $\Vert \nabla u\Vert_{L^p(\B^2)}<\infty$. This gives \eqref{eq:CZ est main thm}.

     For the other direction, if the condition $H(v)\cdot v\neq 0$ for $v\in \S^1$ fails, then \cref{item 1:thm main} fails by \Cref{thm:coercive counterexample}. So, assume $\Tilde{H}(v)\cdot v\neq 0$ for $v\neq 0$ and $\Tilde{H}$ is not linear. Without loss of generality, we may assume $\Tilde{H}(v)\cdot v>0$ for $v\neq 0$ (otherwise, we replace $A$ by $-A$).  
    % By \Cref{prop 3vectors}, there exists $c\in \R$ and nonzero vectors that satisfy \eqref{eq:cond A for A+cJ}. 
    % Again, by continuity of $\Tilde{H}$ we can assume $\Tilde{H}(v)\cdot v>0$ for $v\neq 0$ without loss of generality.
     By \Cref{prop 3vectors}
there are $c\in\R$, non-zero $v_1,v_2,v_3$ and $\theta\in(0,1)$ that satisfy \eqref{eq:cond A for A+cJ}. Set $w_i:=(\Tilde H+cJ)(v_i)$ for $J$ as in \eqref{eq:J}.
Since $Jv\cdot v=0$, $v_i\cdot w_i=\Tilde H(v_i)\cdot v_i>0$. 
% Note that $|(A+cJ)(v)|\le (L+|c|)|v|$ for all $v\in\R^2$, and, by \eqref{eq:limit of A} and the positive homogeneity of $\Tilde H$,
%  $
%  \lim_{t\to\infty}\frac{(A+cJ)(tv_i)}{t}=\Tilde H(v_i)+cJv_i=w_i$ for $ i=1,2,3.$
    % for $J$ as in \eqref{eq:J}. 
    Then, the map $A+cJ:\R^2\to \R^2$ satisfies \Cref{cond:A}.  
    So, we can apply \Cref{th:CZ_failure_A} to find $\bar{q}>1$ such that for any $\Lambda>0$ and $p\in (1,\bar{q})$, there exists a Lipschitz function $u\in W^{1,\infty}_0(\B^2)$ and a map $F\in L^\infty(\B^2, \R^2)$ satisfying \eqref{eq:PDE_thmA} and \eqref{eq:CZ_estimate_thmA}. Since $\div(cJ \nabla u)=0$, we get 
    $$\div(A(\nabla u))=\div (F)\quad \text{in}\;\B^2.$$
    This shows that \cref{item 1:thm main} fails, and this completes the proof.

\end{proof}

%%%%%%%%%%%%%%%%%%%%%%%%%%%%%%%%%%%%%%%%%%
%\newpage
\bibliographystyle{abbrv}
\bibliography{bib}

\end{document}